\documentclass[12pt]{article}

\usepackage{amsmath}
\usepackage{amssymb}
\usepackage{amsthm}

\usepackage{graphicx}
\usepackage{color}
\usepackage{epsfig}
\usepackage[mathscr]{eucal}

\newtheorem{theorem}{Theorem}

\theoremstyle{remark}

\theoremstyle{definition}
\newtheorem{prop}[theorem]{Proposition}
\newtheorem{lemma}[theorem]{Lemma}
\newtheorem{cor}[theorem]{Corollary}

\newtheorem{conjecture}{Conjecture}

\newcommand{\ceil}[1]{\ensuremath{\left \lceil {#1} \right \rceil}}

\newcommand{\cF}{{\cal F}}

\newcommand{\cL}{{\cal L}}

\newcommand{\bv}{{\bf v}}
\newcommand{\bw}{{\bf w}}
\newcommand{\be}{{\bf e}}
\newcommand{\bx}{{\bf x}}
\newcommand{\by}{{\bf y}}

\newcommand{\EH}{Erd\H{o}s-Hajnal}

\DeclareMathOperator{\cis}{cis}
\DeclareMathOperator{\sgn}{sgn}

\DeclareMathOperator{\trans}{trans}
\DeclareMathOperator{\inter}{int}

\allowdisplaybreaks

\title{\EH~Sets and Semigroup Decompositions}
\author{Joshua N. Cooper \thanks{Research supported by NSF Grant DMS-0303272.} \\ \small Courant Institute of Mathematics, NYU \\ \small New York, NY \\ \small \texttt{jcooper@cims.nyu.edu}\\}
\date{\today}

\begin{document}

\maketitle

\begin{abstract} Define a set of lines in $\mathbb{R}^3$ to be ``stacked'' with respect to $\bv \in \mathbb{R}^3$ if, from a vantage point far away in the direction of $\bv$, the lines are linearly ordered by the ``crossing over'' relation.  Given a collection of skew lines and a point $\bv$, we ask, what is the largest stacked subset that must be present among the lines?  This question, which appears in \cite{EHP00}, is intimately related to the well-known \EH~conjecture via the Milnor-Thom theorem.  It was recently resolved by a powerful and very general theorem of Alon, Pach, Pinchasi, Radoi\v{c}i\'{c}, and Sharir (\cite{APPRS04}).  We describe these results and discuss several related issues, including a generalization to ``\EH~sets'' and an intriguing problem concerning the decomposability of semi-algebraic sets: Do all semi-algebraic sets belong to the set algebra generated by semigroups in $\mathbb{R}^d$?  Our main result is a resolution of this question in dimensions 1 and 2.
\end{abstract}

Suppose we have a collection $\cL$ of $n$ lines and a direction $\bv$ in $\mathbb{R}^3$.  Define $T$ to be the set of the pairs $(\ell_1,\ell_2) \in \cL \times \cL$ so that $\ell_1$ ``crosses over'' $\ell_2$ from the perspective of a point very far away in the direction of $\bv$.  If the lines are pairwise skew, then $T$ is a tournament.  In \cite{EHP00}, the authors ask for the size of its largest transitive subtournament: a set of lines which are linearly ordered by the ``crossing over''  relation, i.e., which appear to be ``stacked'' (q.v. Figure 1).

\begin{figure}[h]
\begin{center}
\begin{tabular}{c}
\includegraphics{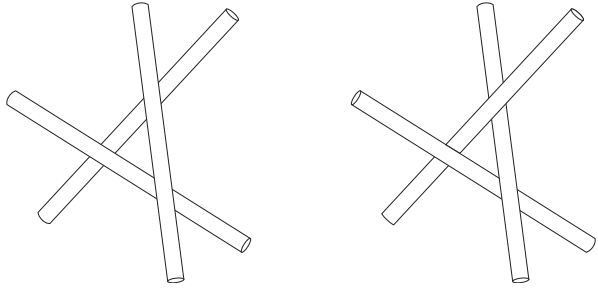}
\end{tabular}
\caption{Stacked lines vs. non-stacked lines.}
\end{center}
\end{figure}

By Ramsey's Theorem, every tournament, including $T$, must have a transitive subtournament of size $\Omega(\log n)$.  Perhaps we can hope for more, however: $T$ is very special, in that it is the result of a very specific construction.  Surely, not all $T$ can arise in this way from line configurations in $\mathbb{R}^3$...?

Indeed, the answer is ``no''.  There are very few such tournaments.  To see this, we first make precise the notion of ``crossing over''.  Parameterize lines in $\mathbb{R}^3$ as follows.  Consider the plane perpendicular to $\bv$ that passes through the origin.  Call this the $xy$-plane.  Then, consider the two planes $x=1$ and $x=-1$.  Each line $\ell$ crosses each plane in one point: $(1,a,b)$ and $(-1,c,d)$, respectively.  (If $\ell$ is parallel to one of them, we may rotate the $xy$-plane and reparameterize so this is no longer the case.  Since $n$ is finite, this can be done for all the lines simultaneously.)  Therefore $\ell$ can be described uniquely by the 4-tuple $(a,b,c,d)$.  Consider two lines, $\ell_1=(a_1,b_1,c_1,d_1)$ and $\ell_2=(a_2,b_2,c_2,d_2)$.  They each project to a line on the $xy$-plane, and cross at a point $(x,y)$, so $\ell_1$ contains a point $(x,y,z_1)$ and $\ell_2$ contains $(x,y,z_2)$.  Whichever $\ell_i$ has the larger $z_i$ ``crosses over'' the other.  We may describe $\ell_i$ as
$$
\{(1,a_i,b_i)+t_i(-2,c_i-a_i,d_i-b_i) : t_i \in \mathbb{R}\},
$$
and its projection $\ell_i^\prime$ onto the $xy$-plane as
$$
\{(1,a_i)+t_i(-2,c_i-a_i) : t_i \in \mathbb{R}\}.
$$
Then $\ell_1^\prime$ and $\ell_2^\prime$ intersect when $t_1 = t_2$ and $t_1(c_1 - a_1)+a_1 = t_1(c_2-a_2)+a_2$, i.e., $t_1 = (a_1-a_2)/((a_1-a_2)-(c_1-c_2))$.  Therefore,
\begin{align*}
z_1 - z_2 & = b_1 - b_2 + (a_1-a_2)\cdot \frac{(d_1-d_2)-(b_1-b_2)}{(a_1-a_2)-(c_1-c_2)}\\
& = \frac{(a_1-a_2)(d_1-d_2)-(b_1-b_2)(c_1-c_2)}{(a_1-a_2)-(c_1-c_2)}.
\end{align*}
Defining $g(a,b,c,d)=(a-c)(ad-bc)$, we have that $\ell_1$ passes over $\ell_2$ iff
$$
g(a_1-a_2,b_1-b_2,c_1-c_2,d_1-d_2) > 0.
$$

To show that this construction can give rise to only a small fraction of all possible tournaments, we employ a standard device: the Milnor-Thom Theorem.  The following version of this important result appears in \cite{W68}.

\begin{theorem}[Milnor '64, Thom '65] Let $R$ be the number of $\pm 1$-vectors in the set 
$$
\left\{\left(g_1(\bx),\ldots,g_m(\bx)\right):\bx \in \mathbb{R}^k\right\}
$$
where $g_i = \sgn \circ f_i$, $1 \leq i \leq m$, and each $f_i$ is a polynomial of degree at most $d$ in $k$ variables.  Then
$$
R \leq \left (\frac{4edm}{k}\right)^k. 
$$
\end{theorem}

In this case, the tournament $T$ is determined by the sign of the $m=\binom{n}{2}$ polynomials $g(a_i-a_j,b_i-b_j,c_i-c_j,d_i-d_j)$, $1 \leq i < j \leq n$, of degree $d=3$ in the $k = 4n$ variables, so the number of such tournaments is at most
$$
\left(3e(n-1)/2\right)^{4n} = 2^{O(n \log n)}.
$$
On the other hand, the total number of tournaments on $n$ vertices is $2^{\binom{n}{2}}$, so this is a vanishingly small fraction.  In particular, there is some $T_0$ which can never be a subtournament of a tournament arising from a line configuration in $\mathbb{R}^3$.

What does this say about the original question, i.e., what is the largest transitive subtournament of a tournament $T$ arising from the crossing patterns of a set of $n$ lines?  Here, the \EH~conjecture plays a role.  Define $\hom(G)$, for a simple graph $G$, to be the size of the largest ``homogeneous subset'' of $G$, i.e., a clique or an independent set.  One way to state Ramsey's Theorem is to say that $\hom(G) \gg \log(n)$ for $G$ on $n$ vertices.

\begin{conjecture}[Erd\H{o}s, Hajnal \cite{EH89}] \label{EHGraphs} For every graph $H$, there is an $\epsilon > 0$ so that, if $G$ on $n$ vertices has no induced copy of $H$, then $\hom(G) \gg n^\epsilon$.
\end{conjecture} 

The intuition behind this conjecture is that the property of avoiding $H$ as an induced subgraph is a mark of being very far from random, since a random graph will have many copies of $H$.  On the other hand, it is notoriously difficult to deterministically construct a graph $G$ coming anywhere close to the Ramsey bound.  In a sense, $\hom(G)$ is a statistic which is {\it very} sensitive to ``true'' randomness.  Therefore, if $G$ avoids $H$, then $\hom(G)$ should be much larger than $\log(n)$ -- something like $n^\epsilon$ for $\epsilon > 0$.

The \EH~conjecture, as stated above, does not directly apply to our situation, however, since we are interested in tournaments instead of graphs.  Fortunately, Alon, Pach, and Solymosi (\cite{APS01}) have provided a very clever argument that the following analogous statement is actually equivalent to Conjecture \ref{EHGraphs}.  Define $\trans(T)$ to be the size of the largest transitive subtournament of the tournament $T$.

\begin{conjecture} \label{EHTournaments} For every tournament $T$, there is an $\epsilon > 0$ so that, if $T^\prime$ on $n$ vertices has no copy of $T$, then $\trans(T^\prime) \gg n^\epsilon$.
\end{conjecture} 

Therefore, if we believe Conjecture \ref{EHGraphs}, then the fact that a tournament arising from a line configuration cannot contain a copy of $T_0$ implies that there is some $\epsilon_\textrm{LINE} > 0$ so that every configuration of lines in $\mathbb{R}^3$ contains a ``stacked'' subset of size $n^{\epsilon_\textrm{LINE}}$.  Indeed, in the same way, the \EH~Conjecture coupled with the Milnor-Thom Theorem implies that {\it any} tournament on a set of points in $\mathbb{R}^d$ whose edges are determined by the sign of a polynomial in the coordinates of the endpoints contains such a ``polynomial-sized'' transitive subtournament.

Recently, a even more general conjecture by Babai in \cite{B78} was resolved unconditionally by Alon, Pach, Pinchasi, Radoi\v{c}i\'{c}, and Sharir (\cite{APPRS04}).  Their powerful result, which has endless applications, is as follows:

\begin{theorem} \label{GEH} Let $\cF$ be a family of $n$ semi-algebraic sets in $\mathbb{R}^d$ of constant description complexity, and let $R \subset \cF \times \cF$ be a fixed semi-algebraic relation on $\cF$.  Then there exists a constant $\epsilon > 0$, which depends only on the maximum description complexity of the sets in $\cF$ and of $R$, and a subfamily $\cF^\prime \subset \cF$ with at least $n^\epsilon$ elements, such that either every pair of distinct elements of $\cF^\prime$ belongs to $R$, or no such pair belongs to $R$.
\end{theorem}

In the case of line configuration tournaments, they are able to show that $\epsilon_\textrm{LINE} \geq 1/6$.  In the other direction, we have:

\begin{prop} $\epsilon_\textrm{LINE} \leq \log 3/\log 7 \approx 0.564575$.
\end{prop}

The proof of this fact relies on a computation and the following observation.  Consider the three ``non-stacked'' lines of Figure 1.  By applying a linear transformation, we may ``bundle'' these lines into the union of a (double) cone and a cylinder, as thin as we wish, without disturbing the crossing relations.  Then, the three cones can be arranged in the same way as the lines of the original three lines, so that now we have $9$ lines with no stacked subset of size larger than $4$.  (See Figure 2.) Proceeding recursively, we can construct a family of $3^k$ lines with no stacked subset of size larger than $2^k$, i.e., $\epsilon_\textrm{LINE} \leq \log 2/\log 3$.

\begin{figure}[h]
\begin{center}
\begin{tabular}{c}
\includegraphics{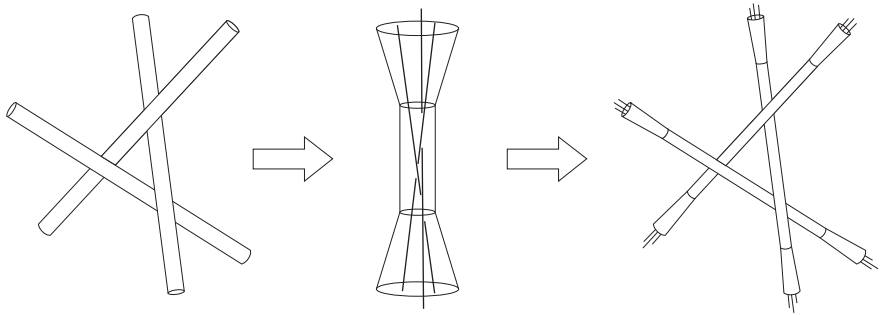}
\end{tabular}
\caption{Bundling lines together recursively.}
\end{center}
\end{figure}

An extensive computer search revealed the existence of a configuration of 7 lines with largest stacked subset of size 3.  Therefore, we can improve this bound by the same construction to $\epsilon_\textrm{LINE} \leq \log 3/\log 7$. \qed \\

Given a real-valued function $f$ on $\mathbb{R}^d$ and a set of points $V \subset \mathbb{R}^d$, define the digraph $\Gamma^f(V)$ to have $V$ as its vertex set and $(x,y) \in \textrm{E}(\Gamma^f(V))$ iff $f(y-x) > 0$.  If $f$ is odd, then $\Gamma^f(S)$ is a tournament; if $f$ is even, we may consider it a simple graph.  By ``complete graph'', we mean the digraph whose edges are all pairs $(u,v)$ with $u$ and $v$ any two distinct vertices.  We call a set $V \subset \mathbb{R}^d$ ``$f$-general'' if $f(v_1-v_2) \neq 0$ for any $v_1, v_2 \in V$.  Call a function $f : \mathbb{R}^d \rightarrow \mathbb{R}$ ``\EH'' if there exists an $\epsilon > 0$ so that, for any $f$-general $V \subset \mathbb{R}^d$ with $|V|=n$, $\Gamma^f(V)$ contains a complete graph, an independent set, or an induced transitive subtournament of cardinality $\gg n^\epsilon$.  The requirement that $V$ be $f$-general corresponds to the requirement that the lines in a configuration be skew.  This assumption, though not strictly necessary, simplifies the discussion below.

It is also possible to generalize this idea as follows.  For a set $S \subset \mathbb{R}^d$ and a set of points $V \subset \mathbb{R}^d$, define $\Gamma^S(V)$ to have $V$ as its vertex set and $(x,y) \in \textrm{E}(\Gamma^S(V))$ iff $y - x \in S$.  Clearly, the definition of $\Gamma^f(V)$ coincides with this definition when $S = f^{-1}(\mathbb{R}^+)$, at least for $f$-general sets $V$.  Furthermore, we say that $S$ is ``\EH'' if there exists an $\epsilon > 0$ so that, for any $V \subset \mathbb{R}^d$ with $|V|=n$, $\Gamma^S(V)$ contains either a complete graph, an independent set, or an induced transitive subtournament of of cardinality $\gg n^\epsilon$.  

Our main question here is: which functions/sets are \EH?  By Theorem \ref{GEH}, all polynomials $f$ and any semi-algebraic set $S$ are \EH.  In the next section, we prove some basic facts and describe a much larger class of \EH~sets.  In particular, we show that all sets belonging to the set algebra $\mathfrak{A}_d$ generated by semigroups in $\mathbb{R}^d$ are \EH.  Therefore, in light of Theorem \ref{GEH}, it is a natural question to ask whether semi-algebraic sets belong to $\mathfrak{A}_n$.  In Section \ref{SectionPolynomial}, we prove that this is the case for $k \leq 2$ and discuss the situation for $k \geq 3$.

\section{\EH~Sets}

We collect some basic results concerning \EH~sets.  First, the set algebra generated by \EH~sets consists entirely of \EH~sets.  For $X \subset \mathbb{R}^d$, write $\bar{X}$ for the topological closure of $X$ and $X^\complement$ for the complement of $X$.

\begin{prop} \label{PropSetAlgebra} Complement, intersection, and union preserve the property of being \EH.
\end{prop}
\begin{proof} The first claim is obvious.  For the second, let $A$ and $B$ be \EH~sets, and consider $G=\Gamma^{A \cap B}(V)$ for some set $V$ in $\mathbb{R}^d$ of cardinality $n$.  The vertex set of $G_1 = \Gamma^{A}(V)$ contains a subset $V_0$ of size $n^{\epsilon_A}$ which induces an independent set, a complete graph, or a transitive tournament.  Since $B$ is \EH, $\Gamma^B(V_0)$ has a subset $V_1$ of size $n^{\epsilon_A \epsilon_B}$ which induces an empty, complete, or transitive digraph $G_2$.  Since $G^\prime = \Gamma^{A \cap B}(V_1) = \Gamma^{A}(V_1) \cap \Gamma^{B}(V_1)$, $G^\prime$ is empty if either $G_1$ or $G_2$ is.  If either $G_i$ is complete, then $G^\prime$ is empty, complete, or transitive according to whether the other $G_j$ is.  If both $G_1$ and $G_2$ are transitive tournaments, then, by the Erd\H{o}s-Szekeres Theorem, their intersection contains either a transitive subtournament of size $n^{\epsilon_A \epsilon_B/2}$ or an independent set of size $n^{\epsilon_A \epsilon_B/2}$.

The fact that the union of two \EH~sets is \EH~now follows from De Morgan's Law.
\end{proof}

\begin{cor} \label{CorMultiply} If $f$ and $g$ are \EH~functions, then $f\cdot g$ is also.
\end{cor}
\begin{proof} The semialgebraic set $(fg)^{-1}(\mathbb{R}^+)$ is the complement of the symmetric difference of the sets $A = f^{-1}(\mathbb{R}^+)$ and $B = g^{-1}(\mathbb{R}^+)$, i.e., $(A \cap B) \cup (A^\complement \cap B^\complement)$, except for some points of $f^{-1}(0)$ and $g^{-1}(0)$.  These sets are of no consequence, however, since $fg$-generality implies $f$-generality and $g$-generality.
\end{proof}

Our next results says that any {\it bounded} set which does not have the origin on its boundary is \EH.

\begin{prop} \label{PropBounded} Any bounded set $B \subset \mathbb{R}^d$ with $\textbf{0} \not \in \partial B$ is \EH.
\end{prop}
\begin{proof} First, suppose that $\textbf{0} \in \bar{B}$.  Since $\textbf{0} \not \in \partial B = \bar{B} \setminus \inter{B}$, it must be that $\textbf{0} \in \inter{B}$.  Therefore, there is some $r > 0$ so that the closed ball of radius $r$ about the origin is contained in $B$, and an $R \geq r$ so that the closed ball of radius $R$ about the origin contains $B$.  Partition space into a lattice of axis-parallel hypercubes (i.e., $[0,s)^d$) with side length $s > 0$ so that $s \sqrt{d} \leq r$.  Let $\Lambda$ denote the set of cubes containing the points of the lattice generated by $t\be_1,\ldots,t\be_d$, where $t = s (\ceil{R/s} + 1)$.  The $t^d$ distinct ``cosets'' $\Lambda_1, \ldots, \Lambda_{t^d}$ of $\Lambda$ partition the points of $V$, and no such $\Lambda_i$ contains a pair $x,y$ with $r < |y-x| \leq R$.  Some $\Lambda_i$ contains at least $n^\prime = n/t^d$ points of $V$, and, for each $x,y \in \Lambda_i$, either there are $\sqrt{n^\prime}$ points $V_0 \subset V$ in one of the cubes comprising $\Lambda_i$, or there are $\sqrt{n^\prime}$ points $V_0 \subset V$ in distinct cubes of $\Lambda_i$.  In the former case, $y-x \in B$ for every pair of points in $V_0$, so $V_0$ induces a complete digraph.  In the latter case, $y-x \not \in B$ for every pair $x,y \in V_0$, so $V_0$ induces an independent set.

On the other hand, if $\textbf{0} \not \in \bar{B}$, then $\textbf{0} \in \inter(B^\complement)$.  There are numbers $r$ and $R$, $0 < r < R$, so that the closed ball $U$ of radius $R$ about the origin contains $B$ and the closed ball $U^\prime$ of radius $r$ about the origin is disjoint from $B$.  Let $A = U \setminus B$, so $U^\prime \subset A$.  $A$ and $U$ are bounded sets whose closure contains $\textbf{0}$, but whose boundary does not.  We may therefore apply the first part of the proof to each of these sets, and invoke Proposition \ref{PropSetAlgebra} by noting that $B = U \setminus A$.
\end{proof}

The last result of this section is that any {\it semigroup} -- a subset of $\mathbb{R}^d$ which is closed under addition -- is \EH.

\begin{prop} \label{PropSemigroup} Any semigroup $S$ is \EH.
\end{prop}
\begin{proof} By Proposition \ref{PropSetAlgebra}, we may assume that $S$ lies in one orthant.  In that case, note that $S \cap -S = \emptyset$.  Furthermore, if $(x,y) \in \Gamma^S(V)$ and $(y,z) \in \Gamma^S(V)$, then $y-x \in S$ and $z-y \in S$, so $z-x = (z-y)+(y-x) \in S$.  Therefore, if we define $y \succ x$ for $x \neq y$ whenever $(x,y) \in \Gamma^S(V)$, then $\succ$ is antisymmetric and transitive, i.e., a partial order.  By Dilworth's Theorem, there is a subset of $V$ of size $n^{1/2}$ which is either a chain or an antichain, i.e., induces a transitive subtournament of $\Gamma^S(V)$ or is an independent set.
\end{proof}

\section{Dimension at Most Two} \label{SectionPolynomial}

Since every element of the set algebra $\mathfrak{A}_d$ generated by semigroups in $\mathbb{R}^d$ is \EH, it is natural to ask whether all semi-algebraic sets belong to this family.  If so, this would give another proof of Theorem \ref{GEH}.  Moreover, we consider this question interesting in its own right.  In this section, we prove the following Theorem.

\begin{theorem} \label{ThmDimensionTwo} Every semi-algebraic set in dimension $d \leq 2$ belongs to $\mathfrak{A}_d$.
\end{theorem}

Note that the statement is obvious in the case that $f$ is a function of only one variable: $S$ is a finite union of intervals whose (possibly infinite) endpoints are on the same side of the origin, and each such interval is either a semigroup itself or a difference of two semigroups.  Also, by Corollary \ref{CorMultiply}, it suffices to prove the theorem for sets defined by irreducible polynomials.  We begin with the following Lemma.

\begin{lemma} \label{LemmaLinear} Any irreducible polynomial $f$ in two variables for which there are infinitely many solutions to $\langle(x,y),\nabla f(x,y)\rangle= 0$ with $f(x,y) = 0$ must be of the form $f(x,y)=ax+by$ for some $a$ and $b$.
\end{lemma}
\begin{proof} Suppose that there are infinitely many pairs $(x,y)$ with $\langle(x,y),\nabla f(x,y)\rangle= 0$ and $f(x,y) = 0$.  By B\'{e}zout's Theorem, the function
$$
h(x,y) = \langle(x,y),\nabla f(x,y)\rangle= x \frac{\partial f}{\partial x} + y \frac{\partial f}{\partial y}
$$
must be divisible by $f(x,y)$.  Consider the terms of highest $x$-degree in $f$, and choose among them the term of highest $y$ degree, say $c x^r y^s$.  It is easy to see that the analogous term for $h$ is $(r+s)cx^ry^s$.  Since these terms are unique, $(r+s)c = c$, i.e., $r+s=1$.  Therefore, $f = ax+g(y)$ for some single-variable polynomial $g$ and constant $a$.  By symmetry, then, $f$ must be of the form $f = ax + by + d$.  However, $f$ cannot have any constant term, since $h$ has no constant term, so $d=0$.
\end{proof}

\begin{lemma} \label{LemmaCircular} Any irreducible polynomial $f$ in two variables for which there are infinitely many solutions to $\langle(-y,x),\nabla f(x,y)\rangle= 0$ with $f(x,y) = 0$ must be of the form $f(x,y)=a(x^2+y^2)+b$ for some $a$ and $b$.
\end{lemma}
\begin{proof} Suppose that there are infinitely many common solutions to $f(x,y)=0$ and $\langle(-y,x),\nabla f(x,y)\rangle=0$.  By B\'{e}zout's Theorem, the function
$$
h(x,y) = \langle(-y,x),\nabla f(x,y)\rangle= x \frac{\partial f}{\partial y} - y \frac{\partial f}{\partial x}
$$
must be divisible by $f(x,y)$.  Write $h = g \cdot f$.  The total degree of $h$ is the same as the total degree of $f$, so $g(x,y) = b$, a constant.  Furthermore, the operator $x D_y - y D_x$ preserves total degree.  Fix a degree $d$, denote by $c_i$ the coefficient of $x^i y^{d-i}$ in $f$, and denote by $c^\prime_i$ the coefficient of $x^i y^{d-i}$ in $h$.  Then
$$
\left [ \begin{array}{c} c^\prime_0 \\ c^\prime_1 \\ \vdots \\ c^\prime_{d-1} \\ c^\prime_d \end{array} \right ] = 
\left [ \begin{array}{ccccccc} 0 & 1 & 0 & \cdots & 0 & 0 & 0 \\ -d & 0 & 2 & \cdots & 0 & 0 & 0 \\ 0 & 1-d & 0 & \ddots & 0 & 0 & 0 \\ \vdots & \vdots & \ddots & \ddots & \ddots & \vdots & \vdots \\ 0 & 0 & 0 & \ddots & 0 & d-1 & 0 \\ 0 & 0 & 0 & \cdots & -2 & 0 & d \\ 0 & 0 & 0 & \cdots & 0 & -1 & 0 \end{array} \right ] \left [ \begin{array}{c} c_0 \\ c_1 \\ \vdots \\ c_{d-1} \\ c_d \end{array} \right ] = 
b \cdot \left [ \begin{array}{c} c_0 \\ c_1 \\ \vdots \\ c_{d-1} \\ c_d \end{array} \right ].
$$
Therefore, $b$ is an eigenvalue of this matrix, which we will denote $M_d$.  We argue that the (nonzero) eigenvalues of $M_d$ are pure imaginary.  Let $D$ be the $d \times d$ diagonal matrix whose $(k+1)^\textrm{th}$ diagonal entry is $i^k \sqrt{\binom{d}{k}}$.  Then $M^\prime_d = D^{-1} M_d D$ is a symmetric tridiagonal matrix whose diagonal entries are all zero and whose $k^\textrm{th}$ super- or subdiagonal element is $i k (-1)^{k+1} \sqrt{\binom{n}{k-1}\binom{n}{k}}$.  Since $M^\prime_d$ is symmetric pure imaginary, its eigenvalues -- and therefore the eigenvalues of $M_d$ -- are pure imaginary.  This contradicts the fact that the $c^\prime_j$ and $c_j$ are all real unless $b=0$.

Now, the last $d$ columns of $M_d$ are linearly independent, so the nullspace of $M_d$ has dimension at most one.  On the other hand, the negative of each eigenvalue $\lambda \in \sigma(M_d)$ is also an eigenvalue, since
$$
M_d (J \bv) = (M_d J) \bv = - (J M_d) \bv = - J (M_d \bv) = - J \lambda \bv =  -\lambda (J \bv),
$$
where $J$ is the ``exchange'' matrix whose only nonzero entries are ones on the antidiagonal, and $\bv$ is an eigenvector corresponding to $\lambda$.  Therefore, the number of nonzero elements of $\sigma(M_d)$ is even.  We may conclude that $M_d$ is invertible if $d$ is odd, and has nullity $1$ if $d$ is even.  In the latter case, it is easy to verify that the $(d+1)$-vector $\bw$ whose $k^\textrm{th}$ element is $0$ if $k$ is even and $\binom{d/2}{(k-1)/2}$ if $k$ is odd is a null vector of $M_d$.  The polynomial corresponding to $\bw$ is $(x^2+y^2)^{d/2}$, so we may write $f = f_0(x^2+y^2)$ for some polynomial $f_0$.

If the factorization of $f_0$ over $\mathbb{R}$ has no linear terms, then $f^{-1}(\mathbb{R}^+) = \emptyset$, a contradiction.  Therefore, $f_0(z)$ is divisible by $z+a$ for some $a \in \mathbb{R}$, and $f$ is divisible by $x^2+y^2+a$.  Since $f$ is irreducible, the lemma follows.
\end{proof}

\begin{lemma} For any irreducible polynomial $f(x,y)$ in two variables which is not of the form $ax+by$ or $a(x^2+y^2)+b$ for any $a,b$, there is a family $\{I_j\}_{j=1}^m$ of open intervals of $[0,2\pi)$, with $|I_j| \leq \pi/2$ for each $j$, and a corresponding family $\{g_j\}_{j=1}^m$ of $C^\infty$, monotone functions defined on each of the $I_j$'s respectively, so that the variety defined by $f$ -- i.e., $\{(x,y):f(x,y)=0\}$ -- is the disjoint union of the graphs of the $g_j$ {\it in polar coordinates}, with the exception of finitely many points.
\end{lemma}

\begin{proof} The set $X$ of zeroes of $f$ is a one-dimensional affine variety with finitely many singularities and finitely many components (in the Euclidean topology).  (See \cite{BCR98} for basic facts concerning the differential topology of real algebraic varieties.)  Furthermore, by the implicit function theorem, removing the singularities and $\textbf{0}$ from $X$ leaves a finite number of $C^\infty$ arcs, i.e., diffeomorphic images of open intervals.  (We consider a topological $S^1$ an arc in the obvious way, i.e., by removing one point.)  Since $f$ is not divisible by $x$ or $y$, $X$ has finitely many points in common with the $x$ and $y$-axes, so we may even remove the two axes and still obtain such a decomposition.  Clearly, it suffices to prove the Lemma for a single one of the resulting arcs, say, $\{(x(t),y(t)):t\in (0,1)\}$.  We may also assume, without loss of generality, that all the points of this arc lie in the first quadrant.

Consider the function $\phi(t) = \tan^{-1}(y(t)/x(t))$.  Differentiating yields
$$
\phi^\prime(t) = \frac{y^{\prime} x - x^\prime y}{y^2 + x^2}.
$$
Suppose there are infinitely many $t \in (0,1)$ for which $\phi^\prime(t) = 0$.  Then $y^{\prime} x - x^\prime y = 0$ at each such $t$, and multiplying by $\partial f/\partial x$ and applying the identity $x^\prime \cdot \partial f/\partial x + y^\prime \cdot \partial f/\partial y = 0$ (obtained by differentiating $f(x,y)=0$),
\begin{align*}
0 & = \frac{\partial f}{\partial x} (y^\prime x - x^\prime y) = \frac{\partial f}{\partial x} y^\prime x - \frac{\partial f}{\partial x} x^\prime y \\
& = \frac{\partial f}{\partial x} y^\prime x + \frac{\partial f}{\partial y} y^\prime y \\
& = y^\prime \cdot \langle (x,y),\nabla f(x,y) \rangle .
\end{align*}
By symmetry, $x^\prime \langle (x,y),\nabla f(x,y) \rangle = 0$.  However, $x^\prime$ and $y^\prime$ cannot both be zero, so $\langle (x,y),\nabla f(x,y) \rangle = 0$.  Lemma \ref{LemmaLinear} implies that $f(x,y)=ax+by$ for some $a$ and $b$, contradicting our hypothesis.  Hence, there exist only a finite set $\{t_j\}_{j=1}^l$ of values of $t \in (0,1)$ for which $\phi^\prime(t)=0$.

Now, consider the quantity $d(x^2+y^2)/d\phi = 2 (xx^\prime + y y^\prime)(x^2+y^2)/(y^\prime x - x^\prime y)$.  If this equals zero at $t \not \in \{t_j\}_{j=1}^l$, then $xx^\prime + y y^\prime = 0$.  Again, applying the identity $x^\prime \cdot \partial f/\partial x + y^\prime \cdot \partial f/\partial y = 0$ and noting that $x^\prime$ and $y^\prime$ cannot simultaneously be zero, we obtain that
$$
x \frac{\partial f}{\partial y} - y  \frac{\partial f}{\partial x} = 0.
$$
By Lemma \ref{LemmaCircular}, there are only finitely many $t^\prime_j$ for which this equality can hold.

If we remove the points $(x(t_j),y(t_j))$ and $(x(t^\prime_j),y(t^\prime_j))$ from the arc, we are left with a finite collection of arcs on which $\phi$ is $C^\infty$, monotone, and defined on an open interval of length at most $\pi/2$.  Choosing one such interval $I$ and letting 
$$
g(\theta) = \sqrt{(x \circ (\phi|_I)^{-1}(\theta))^2 + (y \circ (\phi|_I)^{-1}(\theta))^2}
$$
gives one of the $g_j$ of the desired decomposition: $g_j$ is monotone because
$$
\frac{d g_j^2}{d \theta} = 2 \left((x \circ (\phi|_I)^{-1}(\theta)) \cdot \left . \frac{dx}{d\phi} \right |_{\phi=\theta} + (y \circ (\phi|_I)^{-1}(\theta)) \cdot \left . \frac{dy}{d\phi}\right|_{\phi=\theta}\right) \cdot \frac{d \phi^{-1}}{d \theta} \neq 0.
$$
The other $g_j$ are obtained analogously.
\end{proof}

\begin{proof}[Proof of Theorem \ref{ThmDimensionTwo}] First, if $f(x,y)=ax+by$, then $S$ is a half-plane, which is a semigroup, and the result follows immediately from Proposition \ref{PropSemigroup}.  If $f(x,y) = a(x^2+y^2)+b$, then either $S$ is empty (in which case the result is immediate) or else $S^\complement$ is a disk (in which case $S$ is the union of four semigroups).  Otherwise, let $\{g_j\}_{j=1}^N$ be the collection of all the $g_j$ constructed in the previous Lemma, and let $\{I_j\}_{j=1}^N$ be the corresponding intervals.  Some of the $g_j$ we will classify as ``unimportant'' as follows.  For each point $\bx_t$ on the graph $\Gamma_j$ of $g_j$ -- i.e., the image of $t \in (0,1)$ in polar coordinates -- consider the quantities
$$
T_-(t) = \lim_{\lambda \rightarrow 1-} \sgn(f(\lambda \bx_t)) \, \textrm{ and } \, T_+(t) = \lim_{\lambda \rightarrow 1+} \sgn(f(\lambda \bx_t)).
$$
We argue these limits exist and equal $\pm 1$.  If they did not exist (or equalled 0), then, since $\sgn$ is continuous everywhere but at $0$, there must an infinite sequence of $\lambda_i \rightarrow 1-$ (say) so that $f(\lambda_i \bx_t) = 0$.  But then, by B\'ezout's Theorem, $f$ must be divisible by the equation of the line passing through $\bx_t$ and the origin, a contradiction.  Furthermore, $T_-$ and $T_+$ are continuous functions of $t$ : suppose, to the contrary, that there is a point $\bx_t \in \Gamma_j$ and two sequences $t^1_i$ and $t^2_i$ so that $T_-(t^1_i) = 1$ and $T_-(t^2_i) = -1$, $t^1_i,t^2_i \rightarrow t$.  Then there are sequences $\lambda^1_i$ and a $\lambda^2_i$ so that $\lambda^1_i,\lambda^2_i \rightarrow 1-$, $\sgn(f(\lambda^1_i x^1_i))=+1$ and $\sgn(f(\lambda^2_i x^2_i))=-1$.  For each $i$, we interpolate between $\lambda^1_i x^1_i$ and $\lambda^2_i x^2_i$ with the function $h_i : [t^1_i,t^2_i] \rightarrow \mathbb{R}^2$ defined by $h_i(s) = \lambda_i(s) \bx_s$ where
$$
\lambda_i(s) = \frac{\lambda^2_i - \lambda^1_i}{t^2_i - t^1_i} (s-t^1_i) + \lambda^1_i.
$$
Since $f(h_i(t^1_i)) > 0$ and $f(h_i(t^2_i)) < 0$, there is an $s_i$ so that $f(h_i(s_i)) = 0$, by the intermediate value theorem.  Furthermore, it is easy to see that $s_i \rightarrow t$ and $\lambda_i(s_i) \rightarrow 1-$.  Then there is a sequence of points of $V = f^{-1}(0)$ not lying on the curve $g_j$, but approaching $g_j(t)$, implying that $g_j(t)$ is a singular point of $V$ (a self-intersection, in particular).  This contradicts the fact that $g_j$ contains no singularities of $V$.  Therefore, we may conclude that $T_-$ and $T_+$ are continuous.  Since they take on the values $\pm 1$, however, this means that they are constant.  Hence, we may write $c_j^-$ for the value of $T_-$ on $g_j$ and $c^+_j$ analogously for $T_+$.  Finally, we say that $g_j$ is ``unimportant'' if $c_j^- = c_j^+$, and re-index the $g_j$ so that all the elements of $\{g_j\}_{j=1}^M$ are ``important''.

Now, let $X_j$ denote the set $\{\lambda g_j(\theta) \cis(\theta) : \lambda \geq 1,\, \theta \in I_j\}$ if $c_j^- = +1$ and $\{\lambda g_j(\theta) \cis(\theta) : \lambda > 1,\, \theta \in I_j\}$ if $c_j^- = -1$, where $\cis(\theta) = (\cos(\theta),\sin(\theta))$; let $Y$ be the set of $\bv \in \mathbb{R}^2 \setminus \{\textbf{0}\}$ so that for all sufficiently small $\epsilon > 0$ we have $f(\epsilon \bv) > 0$; let $Y^\prime$ be the set of $\bv \in \mathbb{R}^2 \setminus \{\textbf{0}\}$ so that for all sufficiently small $\epsilon > 0$ we have $f(\epsilon \bv) < 0$; and define $Z$ to be $Y \cdot \mathbb{R}^+$.  Define $A$ to be the set of points lying on the graph of some $g_j$, let $B$ be the (finite) set of points $f^{-1}(0) \setminus A$, let $C = B \cdot \mathbb{R}^+$, and let $D = f^{-1}(\mathbb{R}^+) \cap C$.  We argue that
\begin{equation} \label{EqSymmetricDiff}
S^\prime = S \setminus \{\textbf{0}\} = ((Z \, \circleddash \, X_1 \,\circleddash \cdots \circleddash\, X_M) \setminus C) \cup D,
\end{equation}
where ``$\circleddash$'' denotes symmetric difference (which is associative).  Denote the set $X_1 \circleddash \cdots \circleddash X_M$ by $W$, and consider a ray $\rho_\bv = \bv \cdot \mathbb{R}^+$ with $\bv \not \in C$.  Let $h_\bv(\lambda) = f(\lambda \bv)$ for $\lambda \geq 0$.  Since $f$ is continuous and the graphs of the $g_j$ completely cover the set $f^{-1}(0) \cap C$, $h_\bv$ undergoes a sign change at $\lambda$ iff $\lambda \bv$ lies on the graph of some ``important'' $g_j$.  Therefore, $S^\prime \cap \rho_\bv = \rho_\bv \cap W^\complement$ whenever $\bv \in Y$.  On the other hand, if $\bv \not \in Y$, then $\bv \in Y^\prime$.  If not, then $f$ has zeroes arbitrarily close to the origin along the ray $\rho_\bv$, which we previously argued implies that $f(x,y)=ax+by$ for some $a$ and $b$.  It is clear then that $S^\prime \cap \rho_\bv = W \cap \rho_\bv$ if $\bv \not \in Y$, so that $S^\prime \setminus C = (Z \circleddash W) \setminus C$.

Define $\eta(\theta)$ to be $1$ if $\cis(\theta) \in Y$ and $-1$ if $\cis(\theta) \in Y^\prime$, and let $l(\theta) = \min\{1,\epsilon\}$, where $\epsilon$ is the supremum of the set of reals so that $0 < \epsilon^\prime < \epsilon$ implies
$$
(\sgn \circ f)(\epsilon^\prime \cdot \cis(\theta)) = \eta(\theta).
$$  
Since all the $\epsilon$ are positive and the unit circle is compact, there is a $\delta$ so that $l(\theta) > \delta$ for all $\theta \in [0,2\pi)$.  Therefore, $\eta(\theta) = (\sgn \circ f)(\delta \cdot \cis(\theta))$ for all $\theta$.  If $\eta$ has infinitely many discontinuities, then $f$ has infinitely many zeroes on the circle of radius $\delta$ about the origin, which in turn implies (as previously argued) that $f(x,y)=a(x^2+y^2)+b$ for some $a$ and $b$.  Hence, $\eta^{-1}(1)$ is a finite collection of intervals, and
$$
Z = \{\lambda \cdot \cis(\theta) : \lambda > 0, \, \theta \in \eta^{-1}(1)\}
$$
is a finite union of cones (projective hulls of convex sets), each of which is a semigroup.

Now, suppose $\bx,\by \in X_j$, so that $\bx = \lambda_1 g_j(\theta_1) \cis(\theta)$ and $\by = \lambda_2 g_j(\theta_2) \cis(\theta)$ with $|\theta_1 - \theta_2| \leq \pi/2$, $\theta_1, \theta_2 \in [0,2\pi)$, and $\lambda_1,\lambda_2 \geq 1$.  We may assume without loss of generality that $\bx$ and $\by$ both lie in the first quadrant.  If we write $\bx+\by$ as $(r,\phi)$ in polar coordinates, then is it easy to see that $\theta_1 < \phi < \theta_2$ and $r \geq \max\{\lambda_1 g_j(\theta_1),\lambda_2 g_j(\theta_2)\}$.  By the monotonicity of $g_j$, $r > g_j(\phi)$.  If $\bx+\by \not \in X_j$, then $r \leq g_j(\phi)$, a contradiction.  Hence, $X_j$ is a semigroup.

That $C$ is a finite union of semigroups is obvious.  We can write $D$ as
$$
D = \bigcup_{\bv \in B} (f^{-1}(\mathbb{R}^+) \cap \rho_\bv).
$$
Because $f_0(\lambda) = f(\lambda \bv)$ is a polynomial in a single variable, $f_0$ is \EH.  Furthermore, $(f_0^{-1}(\mathbb{R}^+)\cdot \bv) \cap \rho_\bv =  f^{-1}(\mathbb{R}^+) \cap \rho_\bv$, and $\rho_\bv$ is a semigroup, so $D \in \mathfrak{A}_2$ since $f_0^{-1}(\mathbb{R}^+) \in \mathfrak{A}_1$ implies $f_0^{-1}(\mathbb{R}^+) \cdot \bv \in \mathfrak{A}_2$.

Finally, since each of the terms in (\ref{EqSymmetricDiff}) belongs to $\mathfrak{A}_2$, we may apply Proposition \ref{PropSetAlgebra} to conclude that $S$ is \EH.
\end{proof}

One is tempted to try to generalize Theorem \ref{ThmDimensionTwo} to higher dimensions.  However, we believe this not to be possible, at least for dimension greater than three.  We have the following (perhaps premature) conjecture.

\begin{conjecture} There exist semi-algebraic sets in dimension $d \geq 4$ which do not belong to $\mathfrak{A}_d$, but all semi-algebraic sets defined by three-variable polynomials {\it do} belong to $\mathfrak{A}_3$.
\end{conjecture}

In fact, we believe the function $f(a,b,c,d) = ad-bc$ provides a counterexample for dimension four.  Ironically, showing that $f$ is \EH~is precisely what is needed to show that $\epsilon_\textrm{LINE} > 0$.

\section{Acknowledgements} \label{SectionAck}

Thank you to J\'{a}nos Pach, Ricky Pollack, and Uli Wagner for helpful discussions and ideas.  And a special thank you to J\'{o}zsef Solymosi for his insight and guidance, as well as for introducing me to this wonderful area of research.

\end{document}